\documentclass[twoside,11pt]{article}
\usepackage{graphics,latexsym,amsmath,algorithm,algorithmic,subcaption}
\usepackage{array,url}
\usepackage[usenames]{color}
\usepackage{enumitem}
\usepackage{multirow}
\usepackage{jmlr2e}

\newtheorem{thm}{Theorem}

\newtheorem{lem}{Lemma}

\numberwithin{equation}{section}

\newcommand{\bgamma}{\boldsymbol{\gamma}}
\newcommand{\RR}{\mathbb R}

\newcommand{\by}{\mathbf{y}}

\newcommand{\balpha}{\boldsymbol{\alpha}}
\newcommand {\bfgamma} {\mbox{\boldmath $\gamma$}}

\def\bs{\mathop{\bf s}}
\def\bu{\mathop{\bf u}}

\def\sign{\mathop{\rm sign}}

\def\argmin{\mathop{\rm argmin}}

\ShortHeadings{Kernel-based estimation for PFLM}{Lv and He and Wang}
\firstpageno{1}

\begin{document}

\title{Kernel-based estimation for partially  functional linear model: Minimax rates and  randomized sketches}

\author{\name Shaogao Lv \email lvsg716@swufe.edu.cn\\
\addr College of Statistics and Mathmetics\\
Nanjing Audit University\\
Nanjing, China
\AND
\name Xin He \email he.xin17@mail.shufe.edu.cn\\
\addr School of Statistics and Management\\
Shanghai University of Finance and Economics\\ 
Shanghai, China
       \AND
\name Junhui Wang \email j.h.wang@cityu.edu.hk \\
\addr 	School of Data Science\\
City University of Hong Kong \\
Kowloon Tong, Kowloon, Hong Kong}


\editor{}
\maketitle

\begin{abstract}
This paper considers the partially functional linear model (PFLM) where all predictive features consist of a functional covariate and a high dimensional scalar vector.  Over an infinite dimensional reproducing kernel Hilbert space, the proposed estimation for  PFLM is a least square approach with two mixed regularizations of  a function-norm and an $\ell_1$-norm. Our main task in this paper is to establish the minimax rates for PFLM under high dimensional setting, and the optimal minimax rates of estimation is established by using various techniques in empirical process theory for analyzing kernel classes. In addition, 
we propose an efficient numerical algorithm based on randomized sketches of the kernel matrix.   Several numerical experiments are implemented to support our method and optimization strategy.
\end{abstract}


\section{Introduction}

In the problem of functional linear regression,  a single functional feature $X(\cdot)$ is assumed to be square-integrable  over an interval $\mathcal{T}$, and the classical functional linear regression  between the response $Y$ and $X$ is given as
\begin{align}
Y=\langle X, f^* \rangle_{\mathcal{L}_2}+\varepsilon,
\end{align}
where the inner product $\langle \cdot, \cdot\rangle_{\mathcal{L}_2}$ is defined as $\langle f, g\rangle_{\mathcal{L}_2}:=\int_\mathcal{T}f(t)g(t)dt$ for any $f,g\in \mathcal{L}_2(\mathcal{T})$. Here
$f^*$ is some slope function within $\mathcal{L}_2(\mathcal{T})$ and  $\varepsilon$ denotes a error term with zero-mean. Let ${(Y_i, X_i):\,i = 1,...,n}$ denote independent and identically distributed (i.i.d.) realizations
from the population $(Y, X)$, there is extensive literature on  estimation of the slope function $f^*$, or the value of  $\langle X, f^* \rangle_{\mathcal{L}_2}$.

In practice, it is often the case that a response is affected by both a high-dimensional  scalar vector  and some  random functional variables as predictive features.
These scenarios partially motivates us to study PFLM under  high dimensional setting.
For simplifying the notations,  this paper assumes that $Y$ and $X(\cdot)$ are centered. To be more precise,
we are concerned with partially functional linear regression with the  functional feature $X$  and scalar predictors ${\bf Z}=(Z_1,...,Z_p)^T\in {\cal R}^p$, and a linear model  links the response $Y$ and predictive features ${\bf U}=(X, {\bf Z})$  that
\begin{equation}\label{model}
Y=\langle X, f^* \rangle_{\mathcal{L}_2}+{\bf Z}^T\bgamma^*+\varepsilon,
\end{equation}
where
$\bgamma^*=(\gamma^*_1,...,\gamma^*_p)^T$ denotes the regression coefficients of the scalar covariates, and $\varepsilon$
is a  standard normal variable and independent of $X$ and $\bf Z$. Under the sparse high dimensional setting, a standard assumption is that the cardinality of the active set $S_0:=\{j: \gamma^*_j\neq 0, \ j = 1,...,p\}$ is far less than $p$, while $p$ and $p_0:=|S_0|$  are allowed to diverge as the sample size $n$ increases. In fact, estimation and variable selection issues for 
partially functional linear models have been investigated via FPCA methods by \cite{Shin2012,Lu2014} and \cite{Kong2016}, respectively. 

In this paper,  we focus on a least square regularized estimation for the slope function and the regression coefficients in \eqref{model} under  a kernel-based framework and high dimension setting. The estimators obtained are based on a
combination of  the least-squared loss with a $\ell_1$-type  penalty and the square of a functional norm,  where the former penalty corresponds to the regression coefficients and the latter one is used to control the kernel complexity.  The optimal minimax rates of estimation is established by using various techniques in empirical process theory for analyzing kernel classes, and an efficient numerical algorithm based on randomized sketches of the kernel matrix is implemented to verify our theoretical findings.

\subsection{Our Contributions}
 This paper makes three main contributions to this functional modeling literature.

Our first contribution is  to establish Theorem 1 stating that with high probability, under mild regularity conditions, the prediction error of our procedure under the  squared $L_2$-norm is bounded by $\big(\frac{p_0\log p}{n}+n^{-\frac{2r}{2r+1}}\big)$, where the quantity $r>1/2$ corresponds to the kernel complexity of one composition kernel $K^{1/2}CK^{1/2}$.  The proof of this upper bound involves two different penalties for analyzing  the obtained estimator in high dimensions, and we want to emphasize that it is very hard to prove constraint cone set that has often been used to define a critical condition (constraint eigenvalues constant) for high-dimensional problems \citep*{Bickel2009,Verzelen2012}. To handle this technical difficulty, 
we combine the methods used in \cite{Muller2015} for high dimensional partial linear models  with various techniques in empirical process theory for analyzing kernel classes \citep{Aronszajn1950,Tony2012,Yuan2010,Zhu2014}.

Our second contribution is to establish algorithm-independent minimax lower bounds  under the squared $L_2$ norm. These minimax lower bounds, stated
in Theorem 2, are determined in terms of the metric entropy of the composition kernel $K^{1/2}CK^{1/2}$ and the sparsity structure of high dimensional scalar coefficients.
For the commonly-used kernels, including the Sobolev classes, these lower bounds match our
achievable results, showing optimality of our estimator for PFLM.   It is worthy noting that, the lower
bound of parametric part does not depend on nonparametric smoothness indices,  coinciding with the classical sparse estimation rate in the high dimensional linear models \citep*{Verzelen2012}. By contrast, the lower bound for estimating $f^*$ turns out to be affected by the regression coefficient $\bgamma^*$.
The proof of Theorem 2 is based on characterizing the packing entropies of the class of nonparametric kernel models, interaction between the composition kernel and high dimensional scalar vector, combined with classical information theoretic techniques involving Fano’s inequality and variants \citep*{Yang1999,Van2000,Tsybakov2009}.

Our third contribution is to consider randomized sketches for our original estimation with statistical dimension.  Despite these attractive statistical properties stated as above, the computational complexity of computing our original estimate  prevents it from being routinely used in large-scale problems. In fact,  a standard implementation for any kernel estimator leads to the
time complexity $O(n^3)$ and space complexity  $O(n^2)$ respectively. To this end, we employ
the random projection and sketching techniques developed in  \cite{Yang2017,Mahoney2011}, where it is proposed to approximate $n$-dimensional kernel matrix by projecting its row and column subspaces to a randomly chosen m-dimensional subspace with $m\ll n$. We give the sketch  dimension $m$ proportional to the statistical dimension, under which  the resulting estimator has a comparable numerical performance.

\subsection{Related Work}
A class of conventional estimation procedures for functional linear regressions in the statistical literature are based on functional principal components regression (FPCA) or spline functions; see \citep*{Ramsay2005,Ferraty2006,Kong2016} and \citep*{Cardot2003} for details. These truncation approaches to handle an infinity-dimensional function only depend on the information of the feature $X$. In particular, commonly-used FPCA methods that form an available basis for the slope function $f^*$ are determined solely by empirical covariance of the observed feature $X$, and these basis may not act as an efficient representation to approximate $f^*$, since  the slope function $f^*$ and the leading functional  components are essentially unrelated. Similar problems also arise when spline-based finite representation are used.

To avoid inappropriate representation for the slope function, reproducing kernel  methods have been known to be a family of powerful tools for directly estimating infinity-dimensional functions. When the slope function is assumed to reside in a reproducing kernel Hilbert Space (RKHS), denoted by $\mathcal{H}_{K}$, several existing work \citep*{Yuan2010,Tony2012,Zhu2014} for functional linear or additive regression have proved that the minimax rate of convergence depends on both the kernel $K$ and the covariance function $C$ of the functional feature $X$. In particular, the
alignment of $K$ and $C$ can significantly affect the optimal rate of convergence. However, it is well known that  kernel-based methods suffer a lot from storage cost and computational burden. Specially, kernel-based methods need to store a $n\times n$ matrix before running algorithms and  are limited to small-scale problems.


\subsection{Paper organization}
The rest of this paper is organized as follows.  Section \ref{sec:2}  introduces some notations and the basic knowledge on
kernel method, and formulates the proposed kernel-based regularized estimation method. Section \ref{sec:3} is devoted to establish the minimax rate of the prediction problem for PFLM and provide detailed discussion on the obtained results, including  the desired convergence rate of the upper bounds and a matching set of
minimax lower bounds. In Section 4,  a general sketching-based strategy is provided, and an approximate algorithm for solving \eqref{matrixopt} is employed.  Several numerical experiments are implemented in Section \ref{sec:num} to support the proposed approach  and the employed optimization strategy. A brief summary of this paper is provided in Section \ref{sec:con}. Appendix A contains several core proof procedures  of the main results, including the technical proofs of Theorems \ref{nonpara}--\ref{thm3}. Some useful lemmas and more technical details  are provided in Appendix B.

\section{ Problem Statement and Proposed Method}\label{sec:2}

\subsection{Notation}
Let $u,v$ be two general random variables, and  denote the joint distribution of $(u,v)$ by $Q$ and  the marginal distribution of $u(z)$ by $Q_u(Q_v)$.
For a measurable function $f:\,u\times v\rightarrow \RR$, we define the squared $L_2$-norm by $\|f\|^2:=\mathbb{E}_Qf^2(u,v)$, and the squared empirical norm is given by $\|f\|_n^2:=\frac{1}{n}\sum_{i=1}^nf^2(u_i,v_i)$, where $\{(u_i,v_i)\}_{i=1}^n$ are i.i.d. copies of $(u,v)$. Note that $Q$ may differ from line to line.
For a vector $\bgamma\in \RR^p$, the $\ell_1$-norm and $\ell_2$-norm are given by $\|\bgamma\|_1:=\sum_{j=1}^p|\gamma_j|$ and $\|\bgamma\|_2:=\big(\sum_{j=1}^p\gamma_j^2\big)^{1/2}$, respectively. With a slight abuse of notation, we write $\|f\|_{\mathcal{L}_2}:=\langle f, f\rangle_{\mathcal{L}_2}$ with $\langle f, g\rangle_{\mathcal{L}_2}=\int_\mathcal{T}f(t)g(t)dt$. For two sequences $\{a_k: k\geq 1\}$ and $\{b_k: k\geq 1\}$, $a_k\lesssim b_k$ (or $a_k=O(b_k)$) means that there exists some constant $c$ such that $a_k\leq cb_k$ for all $k\geq 1$. 
Also, we write $a_k\gtrsim b_k$ if there is some positive constant $c$  such that $ a_k\geq c b_k$ for all $k\geq 1$.
Accordingly, we write $a_k\asymp b_k$ if both $a_k\lesssim b_k$ and $a_k\gtrsim b_k$  are satisfied.

\subsection{Kernel Method} 
Kernel methods are one of the most powerful learning schemes in machine learning, which often take the form of  regularization
schemes in a reproducing kernel Hilbert space (RKHS) associated with a Mercer kernel \citep*{Aronszajn1950}. 
A major advantage of employing the kernel methods is that the corresponding optimization task over an infinite dimensional RKHS are equivalent to a $n$-dimensional optimization problems, benefiting from the so-called reproducing property.

Recall that
a kernel $K(\cdot,\cdot):\mathcal{T} \times \mathcal{T}\rightarrow {\cal R}$ is a continuous, symmetric, and positive semi-definite function.
Let $\mathcal{H}_{{K}}$  be the closure of the linear span of functions $ \{K_t (\cdot):= {K}(t, \cdot), t \in \mathcal{T}  \}$
endowed with the  inner product $\langle \sum_{i=1}^n \alpha_i{K}_{t_i},\,
\sum_{j=1}^n \beta_j{K}_{t_j}\rangle_{{K}}:=
\sum_{i,j=1}^n\alpha_i\beta_j{K}(t_i,t_j)$, for any $\{t_i\}_{i=1}^n, \{t_i\}_{i=1}^n \in \mathcal{T}^n$ and $n\in \mathcal{N^+}$.
An important  property on $\mathcal{H}_{{K}}$ is the reproducing property stating that
$$
f(t)=\langle f, {K}_t\rangle_{{K}},\,\,\, \hbox{for any}\,f\in \mathcal{H}_{{K}}.
$$
This  property ensures that an RKHS inherits many nice properties  from the standard finite dimensional Euclidean spaces. Throughout this paper, we assume that the slope function $f^*$ resides in a specified RKHS, still denoted by 
$\mathcal{H}_K$. In addition, another RKHS can be naturally induced by the stochastic process of $X(\cdot)$.
Without loss of generality, we assume that $X(\cdot)$  is square integrable over $\mathcal{T}$ with zero-mean, ant thus the covariance function of $X$, defining as
$$
C(s,t)=\mathbb{E}[X(s)X(t)],\quad \forall \, t,\,s\in \mathcal{T},
$$
is also a real, semi-definite kernel.

Note that the kernel complexity is characterized explicitly by a kernel-induced integral operator. Precisely, 
for any  kernel ${K(\cdot,\cdot)}: \mathcal{T}\times \mathcal{T}\rightarrow {\cal R}$, we define the integral operator $L_{{K}}: \mathcal{L}_2(\mathcal{T})\rightarrow \mathcal{L}_2(\mathcal{T})$ 
by
$$
L_{{K}}(f)(\cdot)=\int_\mathcal{T} {K}(s,\cdot) f(s)ds.
$$
By the reproducing property, $L_{{K}}$ can be equivalently defined as
$$
\langle f, L_{{K}}(g)\rangle_{K}=\langle f, g\rangle_{\mathcal{L}_2},\quad \forall\,
f\in \mathcal{H}_{{K}},\,g\in\mathcal{L}_2(\mathcal{T}). 
$$
Since the operator $L_{{K}}$ is linear, bounded and self-adjoint in $\mathcal{L}_2(\mathcal{T})$, 
the spectral theorem implies that there exists a family of orthonormalized eigenfunctions $\{\phi^{{K}}_\ell:\,\ell\geq 1\}$ and
a sequence of eigenvalues $\theta_1^{{K}}\geq \theta_2^{{K}}\geq ...>0$ such that
$$
{K}(s,t)=\sum_{\ell\geq 1} \theta_\ell^{{K}} \phi^{{K}}_\ell(s)\phi^{{K}}_\ell(t),\quad s,\,t\in \mathcal{T},
$$
and thus by definition, it holds
$$
L_{{K}}(\phi^{{K}}_\ell)= \theta_\ell^{{K}}\phi^{{K}}_\ell,\quad \ell=1,2,...
$$
Based on the semi-definiteness of  $L_{{K}}$, we can always decompose it into the following form
$$
L_{{K}}=L_{{K}}^{1/2}\circ L_{{K}}^{1/2},
$$
where  $L_{{K}^{1/2}}$ is also a  kernel-induced  integral operator associated with a fractional kernel ${K}^{1/2}$  that
$$
{K}^{1/2}(s,t):=\sum_{\ell\geq 1} \sqrt{\theta_\ell^{{K}}} \phi^{{K}}_\ell(s)\phi^{{K}}_\ell(t),\quad s,\,t\in \mathcal{T}.
$$
Also, it holds
$$
L_{{K}^{1/2}}(\phi^{K}_\ell):=\sqrt{\theta_\ell^{{K}}}\phi^{{K}}_\ell.
$$
Given two kernels $K_1, K_2$, we define
$$
(K_1K_2)(s,t):=\int_\mathcal{T} K_1(s,u) K_2(t,u)du,
$$
and then it holds $L_{K_1K_2}=L_{K_1}\circ L_{K_2}$. Note that $K_1K_2$ is not necessarily a symmetric  kernel.

In the rest of this paper, we focus on the RKHS $\mathcal{H}_{K}$ in which the slope function $f^*$ in \eqref{model} resides. 
Given the kernel $K$, the covariance function $C$ and by using the above notation, we define the linear operator $L_{K^{1/2}C_kK^{1/2}}$  by 
$$
L_{K^{1/2}CK^{1/2}}:=L_{K^{1/2}}\circ L_{C}\circ L_{K^{1/2}}.
$$
If the both  operators $L_{K^{1/2}}$ and $L_{C}$ are  linear, bounded and self-adjoint, so is $L_{K^{1/2}CK^{1/2}}$. By the spectral theorem, there exist a sequence of
positive eigenvalues $s_{1}\geq s_{2}\geq ...>0$ and a set of orthonormalied eigenfunctions $\{\varphi_{\ell}: \ell\geq 1\}$ such that
$$
K^{1/2}CK^{1/2}(s,t)=\sum_{\ell\geq 1}s_{\ell} \varphi_{\ell}(s)\varphi_{\ell}(t),\quad \forall\, s,t \in \mathcal{T},
$$
and particularly
$$
L_{K^{1/2}CK^{1/2}}(\varphi_\ell)=s_{\ell} \varphi_{\ell}, \quad \ell=1,2,...
$$
It is worthwhile to note that the eigenvalues $\{s_{\ell}: \ell\geq 1\}$ of the linear operator $L_{K^{1/2}CK^{1/2}}$ depend on the eigenvalues of both
the reproducing kernel $K$ and the covariance function $C$. We shall show in Section \ref{sec:3} that the minimax rate of convergence of the excess prediction risk
is determined by the decay rate of the eigenvalues $\{s_{\ell}: \ell\geq 1\}$.

\subsection{Regularized Estimation and Randomized Sketches}

Given the sample $\{Y_i,(X_i, {\bf Z}_i)\}_{i=1}^n$ which are independently drawn from \eqref{model},  the proposed estimation procedure for PFLM \eqref{model} is formulated in a least square regularization scheme by solving
\begin{equation}\label{method}
(\widehat{f}, \widehat{\bgamma})=\argmin_{f\in \mathcal{H}_K,\bgamma \in {\cal R}^p}\Big\{\frac{1}{n}\sum_{i=1}^n\big(Y_i-\langle X_{i},f \rangle_{\mathcal{L}_2} -{\bf Z}^T_i\bgamma\big)^2+\mu^2\|f\|^2_{K}+\lambda \|\bgamma\|_1\Big\},
\end{equation}
where  the parameter $\mu^2 > 0$ is used to control the smoothness  of
the nonparametric component and  $\lambda>0$ associated with the $\ell_1$-type penalty  is used to generate sparsity with respect to the scalar covariates.

Note that although the proposed estimation procedure \eqref{method} is  formulated  within an infinite-dimensional Hilbert space, the following lemma shows that this optimization task
is equivalent to a finite-dimensional minimization problem.
\begin{lem}\label{finiteexpre}
	The proposed estimation procedure  \eqref{method} defined on $\mathcal{H}_{K} \times {\cal R}^p$ is equivalent to a finite-dimensional parametric convex optimization. That is, $\widehat{f}(t)=\sum_{k=1}^n \alpha_kB_k(t)$ with unknown coefficients $\balpha=(\alpha_1,...,\alpha_n)^T$, for any $t\in \mathcal{T}$. Here each basis function $B_k(t)=\langle X_k, K(t,) \rangle_{\mathcal{L}_2(\mathcal{T})}\in \mathcal{H}_{K}$, $k=1,...,n$.
\end{lem}

To rewrite the minimization problem \eqref{method} into a matrix form, we define a $n\times n$ semi-definite matrix $\mathbb{K}^c=(K^c_{ik})_{i,k=1}^n$ with $K^c_{ik}:=\langle X_{i},B_k \rangle_{\mathcal{L}_2}=\iint X_k(u)X_i(t)K(t,u)dudt$, and by the reproducing property on $K$, we also get
$\langle B_{i},B_k \rangle_{K}=K^c_{ik}$, $i,k=1,...,n$.  Thus, by Lemma \ref{finiteexpre},  the matrix form of \eqref{method} is given as
\begin{align}\label{matrixopt}
\min_{\balpha\in {\cal R}^n,\bgamma\in {\cal R}^p}\frac{1}{n}\big\|\by-\mathbb{K}^c\balpha-\mathbb{Z}\bgamma\big\|_2^2+\mu^2\balpha^T\mathbb{K}^c\balpha+\lambda\|\bgamma\|_1,
\end{align}
where $\mathbb{Z}\in {\cal R}^{n\times p}$ denotes the design matrix of $\bf Z$.
Since the unconstrained problem \eqref{matrixopt} is convex for both $\balpha$ and $\bgamma$, the standard alternative optimization \citep*{Boyd2004} can be applied directly to approximate a global minimizer of \eqref{matrixopt}. Yet, due to the fact that $\mathbb{K}^c$ is a $n\times n$  matrix, both  computation cost and  storage burden are very heavy in standard implementation, with the orders $O(n^3)$ and $O(n^2)$ respectively. To alleviate the computational issue, we propose an approximate numerical optimization instead of \eqref{matrixopt} in Section 4. Precisely, a class of general random projections are adopted to compress the original kernel matrix $\mathbb{K}^c$ and improve the computational efficiency. 

\section{Main Results: Minimax Rates}\label{sec:3}

In this section, we present the main theoretical results of the proposed estimation in the minimax sense. Specifically, we derive the minimax rates  in term of prediction error  for the estimators in \eqref{method} under high dimensional and kernel-based frameworks. The first two theorems prove the convergence of the obtained estimators, while the last one provides an algorithmic-independent lower bound for the prediction error. 

\subsection{Upper Bounds}
We denote the short-hand notation 
$$
\mathcal{G}:=\big\{g=\langle X,f \rangle_{\mathcal{L}_2}+{\bf Z}^T\bgamma,\,\,f\in \mathcal{H}_{K},\,\bgamma\in {\cal R}^p\big\},
$$
and the functional $g^*({\bf U}):=\langle X,f^* \rangle_{\mathcal{L}_2}+{\bf Z}^T\bgamma^*$ for ${\bf U}=(X,{\bf Z})$. 
With a slight confusion of notation,  we sometimes also write $\mathcal{G}:=\big\{g=(f,\bgamma),\,\,f\in \mathcal{H}_{K},\,\bgamma\in {\cal R}^p\big\}$. 
To split the scalar components and the functional component involved in our analysis, we define the projection of $Z_j$ concerning $\mathcal{H}_{K}$ as
$\Pi(Z_j|\mathcal{H}_{K})=\argmin_{f\in \mathcal{H}_{K}}\|Z_j-\langle X,f \rangle_{\mathcal{L}_2}\|^2$. Let $\Pi(Z_j|X)=\langle X,\Pi(Z_j|\mathcal{H}_{K}) \rangle_{\mathcal{L}_2}$ and  $\Pi_{{\bf Z}|X}=(\Pi(Z_1|X),...,\Pi(Z_p|X))^T$, and then we denote  $\widetilde{\bf Z}:={\bf Z}-\Pi_{{\bf Z}|X}$ as a random vector of ${\cal R}^p$. 
For any $g_1({\bf U}):=\langle X,f_1 \rangle_{\mathcal{L}_2}+{\bf Z}^T\bgamma_1\in \mathcal{G}$ and $g_2({\bf U}):=\langle X,f_2 \rangle_{\mathcal{L}_2}+{\bf Z}^T\bgamma_2\in \mathcal{G}$, we have the following orthogonal decomposition that
\begin{align*}
g_1({\bf U})-g_2({\bf U})&={\bf Z}^T(\bgamma_1-\bgamma_2)+\langle X,f_2-f_1 \rangle_{\mathcal{L}_2} \\
&=\widetilde{\bf Z}^T(\bgamma_1-\bgamma_2)+\Pi_{{\bf Z}^T|X}^T(\bgamma_1-\bgamma_2)+\langle X,f_2-f_1 \rangle_{\mathcal{L}_2},
\end{align*}
and by the definition of projection, it holds
\begin{align}\label{othdecom}
\|g_1-g_2\|^2=\|\widetilde{\bf Z}^T(\bgamma_1-\bgamma_2)\|^2+\|\Pi_{{\bf Z}|X}^T(\bgamma_1-\bgamma_2)+\langle X,f_2-f_1 \rangle_{\mathcal{L}_2}\|^2.
\end{align}

To establish the refined  upper bounds of the prediction  and estimation errors, we summarize and discuss the main conditions needed in the theoretical analysis below.

\noindent{\bf Condition A}(Eigenvalues condition). The smallest eigenvalue $\Lambda^2_{min}$ of {$\mathbb{E}[\widetilde{\bf Z}\widetilde{\bf Z}^T]$} is positive,
and the largest eigenvalue $\Lambda^2_{max}$ of {$\mathbb{E}[\Pi_{{\bf Z}|X}\Pi_{{\bf Z}|X}^T]$} is finite.

\noindent {\bf Condition B}(Design condition). For some positive constants $C_z,C_\pi,C_h$, there  holds:
$$
|Z_{j}|\leq C_z, \,\, \|\Pi(Z_j|X)\|_\infty\leq C_\pi,\,
\hbox{and}\,\, \|\Pi(Z_j|\mathcal{H}_{K})\|_{K}\leq C_h,\quad \mbox{for any}  \,j=1,...,p.
$$

\noindent {\bf Condition C}(Light tail condition). There exist two constants $c_1,\,c_2$ such that $$
\mathbb{P}\{\|L_{K^{1/2}}X\|_{\mathcal{L}_2}\geq t\}
\leq c_1\exp(-c_2t^2),\quad  \mbox{for any} \ t>0.
$$


\noindent {\bf Condition D}(Entropy condition). For some constant $1/2<r<\infty$, 
the sequence of eigenvalues $s_\ell$ satisfy that
$$
s_\ell \asymp \ell^{-2r},\quad \ell \in \mathbb{N}^+.
$$

Condition A is commonly used in  literature of semiparametric modelling ; see \citep*{Muller2015} for reference.
This condition ensures that there is enough information in the data to identify the parameters
in the scalar part. Condition B imposes some boundedness assumptions, which are not essential and are used only for simplifying the technical proofs.  Condition C implies that the random process $L_{K^{1/2}}X$ has an exponential decay rate and the same condition  is also considered in \cite{Tony2012}. Particularly,  it is naturally satisfied  if $X$ is a Gaussian process. 
In Condition D,  the parameters $s_\ell$ are related to the alignment between $K$ and $C$, which plays an important role in determining the minimax optimal rates. Moreover, the decay of $s_\ell$ characterizes the kernel complexity and has close relation with various covering numbers and Radmeacher complexity.  Specially, the polynomial decay  assumed in  Condition D  is satisfied for the classical Sobolev class and Besov class.

The following theorem states that with an appropriately chosen $(\mu,\lambda)$, the  predictor $\widehat{g}:=\langle X,\widehat{f}\rangle_{\mathcal{L}_2}+{\bf Z}^T\widehat{\bgamma}$ attains a sharp convergence rate {under  $L_2$-norm}. 

\begin{thm}\label{nonpara}
	Suppose that Conditions A-D hold. With the choice of the tuning parameters $(\mu,\lambda)$, such that
	$$
	\mu\asymp n^{-\frac{r}{2r+1}}+\sqrt{\log(2p)/n},\,\, \lambda\asymp\sqrt{\log(2p)/n}.
	$$
	Then with probability at least $1-2\exp[-n(\delta_1'')^2\mu^2]$, the proposed estimation for PFLM satisfies
	$$
	\|\widehat g-g^*\|^2\lesssim\Big( n^{-\frac{2r}{2r+1}}+\frac{p_0\log (2p)}{n}\Big),
	$$
	where the constant $\delta''$ is  small appropriately. 
\end{thm}
Theorem \ref{nonpara} shows that the proposed estimation \eqref{method} achieves a fast convergence rate in the term of  prediction error. Note that the derived rate depends on the kernel complexity  of $K^{1/2} CK^{1/2}$ and the sparsity of  scalar components. It is  interesting to note that even there exists some underlying correlation structure  between the functional feature and the scalar covariates,   the choice of  hyper-parameter $\mu$ depends on  the structural information of all the features, while the sparsity hyper-parameter $\lambda$ only depends on the scalar component.

\begin{thm}\label{upperb}
	Suppose that all the conditions in Theorem \ref{nonpara} are satisfied.
	Then with probability at least $1-4\exp[-n(\delta_1'')^2\mu^2]-\frac{5}{2p}$, there holds
	\begin{eqnarray}\label{kerydet}
	\|\widetilde{\bf Z}^T(\widehat \bgamma-\bgamma^*)\|^2+
	\frac{\lambda}{8}\|\widehat \bgamma-\bgamma^*\|_1
	\lesssim\Big(\frac{p_0}{\Lambda^2_{min}}\frac{\log (2p)}{n}\Big),
	\end{eqnarray}
	and in addition, we have
	\begin{align}
	\|\langle X,\widehat g-g^* \rangle_{\mathcal{L}_2}\|^2\lesssim\Big( n^{-\frac{2r}{2r+1}}+\frac{p_0\log (2p)}{n}\Big).
	\end{align}
\end{thm}

It is worthy pointing out that the estimation error of the parametric estimator $\widehat \bgamma$ can achieve the optimal convergence rate in the high dimensional linear models \citep*{Verzelen2012}, even in the presence of   nonparametric components. {This result in the functional literature is similar in spirit to the classical  high dimensional partial linear models \citep*{Muller2015,Yu2018}.}

\subsection{Lower Bounds}
In this part, we establish the lower bounds on the minimax risk of estimating $\bgamma^*$ and $\langle X,f^*\rangle_{\mathcal{L}_2}$ separately. Let $B[p_0,p]$ be a set of $p$-dimensional vectors with at most $p_0$ non-zero coordinates, and  ${\cal B}_K$ be the unit ball of $\mathcal{H}_{K}$.
Moreover, we define the risk of estimating $\bgamma^*$ as
$$
R_{\bgamma^*}(p_0,p, {\cal B}_K):=\inf_{\widehat{\bgamma}}\sup_{\bgamma^*\in B[p_0,p],f^*\in {\cal B}_K}\mathbb{E}[\|\widehat \bgamma-\bgamma^*\|_2^2],
$$ 
where ${\inf}$ is taken over all possible estimators for $\bgamma^*$ in model \eqref{model}. Similarly, we define
the risk of estimating $\langle X,f^*\rangle_{\mathcal{L}_2}$ as 
$$
R_{f^*}(s_0,p, {\cal B}_K):=\inf_{\widehat{f}}\sup_{\bgamma^*\in B[p_0,p],f^*\in {\cal B}_K}\mathbb{E}[\langle X,\widehat{f}-f^* \rangle_{\mathcal{L}_2}^2]=\inf_{\widehat{f}}\sup_{\bgamma^*\in B[p_0,p],f^*\in {\cal B}_K}\|L_{C^{1/2}}(\widehat{f}-f^*)\|_{\mathcal{L}_2}^2.
$$ 

The following theorem provides the lower bounds of the minimax optimal estimation error for $\bgamma^*$ and the predictor error for $f^*$, respectively. 
\begin{thm}\label{thm3}
	Given $n$ i.i.d. samples from \eqref{model} with the entropy condition (Condition D). When $p$ is diverging as $n$ increases and $p_0\ll p$, the minimax risk for estimating $\bgamma^*$ can be bounded from below as
	$$
	R_{\bgamma^*}(p_0,p, {\cal B}_K)\gtrsim \frac{p_0\log (p/p_0)}{n};
	$$ 
	the minimax  risk for estimating  $\langle X,f^*\rangle_{\mathcal{L}_2}$ can be bounded from below as 
	$$
	R_{f^*}(p_0,p, {\cal B}_K)\gtrsim \max\Big\{\frac{p_0\log (p/p_0)}{n}, n^{-\frac{2r}{2r+1}}\Big\}.
	$$
\end{thm}
The proof of  Theorem \ref{thm3} is provided in Appendix A. As mentioned previously, these  results indicate that the best possible estimation of $\bgamma^*$ is not affected by the existence of nonparametric components, while the
minimax risk for estimating the (nonparametric) slope function not only depends on the smoothness itself, but also on the dimensionality and sparsity of the scalar covariates. From the lower bound of $R_{f^*}(p_0,p,{\cal B}_K)$, we observe that  a rate-switching phenomenon occurring between a sparse regime and a smooth regime. 
Particularly when $\frac{p_0\log (p/p_0)}{n}$ dominates $n^{-\frac{2r}{2r+1}}$ corresponding to 
the sparse regime,  the lower bound becomes the classical high dimensional parametric rate $\frac{p_0\log (p/p_0)}{n}$.
Otherwise, this corresponds to the smooth regime and thus has similar behaviors as classical nonparametric models. We also notice that the minimax lower bound obtained for the predictor error generalizes the previous results for the pure functional linar model \citep*{Tony2012}.

\section{Randomized Sketches and Optimization}\label{sec:alg}

This section is devoted to  considering an approximate algorithm for \eqref{matrixopt}, based on constraining the original parameter $\balpha\in {\cal R}^n$ to an $m$-dimensional subspace of ${\cal R}^n$, where $m \ll n $ is the projection dimension. We define this approximation via a sketch matrix $\mathbb{S} \in {\cal R}^{m\times n}$ such that the $m$-dimensional subspace is generated by the row span of $\mathbb{S}$. More precisely, the sketched kernel partial functional estimator is given by first solving 
\begin{align}\label{sketchopt}
(\widehat{\balpha}_s,\widehat{\bgamma}_s):&=\arg\min_{\balpha\in {\cal R}^m,\bgamma\in {\cal R}^p}\frac{1}{n}\balpha(\mathbb{S}\mathbb{K}^c)(\mathbb{S}\mathbb{K}^c)^T\balpha-
\frac{2}{n}\balpha^T\mathbb{S}\mathbb{K}^c(\by-\mathbb{Z}\bgamma)
+\frac{1}{n}\|\by-\mathbb{Z}\bgamma\|_2^2\nonumber\\
&+\mu^2\balpha^T\mathbb{S}\mathbb{K}^c\mathbb{S}^T\balpha+\lambda\|\bgamma\|_1.
\end{align}
Then the resulting predictor for the slope function $f^*$ is given as
$$
\widehat{f}_s(t):=\sum_{k=1}^n (\mathbb{S}^T\widehat{\balpha}_s)_kB_k(t)=     \widehat{\balpha}_s^T\mathbb{S}{\bf B}(t), \quad \forall\, t\in \mathcal{T}.
$$
where ${\bf B}(t)=(B_1(t),...,B_n(t))^T\in {\cal R}^n$, where $B_k(t)$ is defined in Lemma \ref{finiteexpre}.  By doing randomized sketches, an approximate form of the kernel estimate $\widehat{\balpha}_s$ can be obtained by solving an $m$-dimensional quadratic program when   $\widehat{\bgamma}_s$ is fixed, which involves time and space complexity $O(m^3)$ and $O(m^2
)$. Computing the approximate kernel matrix is a preprocessing step with time complexity $O(n^2\log(m))$ for properly chosen projections.

\subsection{Alternating Optimization}

This section provides the detailed computational issues of the proposed approach. Precisely, we aim to solve the following  optimization task that
\begin{align}\label{comp:1}
&(\widehat{\balpha}_s,\widehat{\bgamma}_s):=\argmin_{\balpha\in {\cal R}^m,\bgamma\in {\cal R}^p}\frac{1}{n}\balpha^T(\mathbb{ S}\mathbb{K}^c)(\mathbb{ S}\mathbb{K}^c)^T\balpha  {-
	\frac{2}{n}\balpha^T\mathbb{ S}  \mathbb{K}^c (\by-\mathbb{Z}\bgamma)}+\nonumber\\
&\hspace{4cm}  {\frac{1}{n}(\by-\mathbb{Z}\bgamma)^T(\by-\mathbb{Z}\bgamma)}  + \mu^2\balpha^T\mathbb{ S}\mathbb{K}^c\mathbb{ S}^T\balpha    +  \lambda\|\bgamma\|_1.
\end{align}
To solve (\ref{comp:1}), a splitting algorithm with proximal operator is applied, which updates  the representer coefficients ${\balpha}$ and the linear coefficients ${\bgamma}$ sequentially. Specifically, at the $t$-th iteration with current solution $(\balpha^t, \bfgamma^t)$, the following two optimization tasks are solved sequentially to obtain  the solution of the $(t+1)$-th iteration 
\begin{align}
&\balpha^{t+1}=\argmin_{\balpha\in {\cal R}^m}  \Big \{ \frac{1}{n}\balpha^T(\mathbb{ S}\mathbb{K}^c)(\mathbb{ S}\mathbb{K}^c)^T\balpha-
\frac{2}{n}\balpha^T\mathbb{ S}   \mathbb{K}^c   (\by-\mathbb{Z}\bgamma^t)          +\mu^2\balpha^T\mathbb{ S}\mathbb{K}^c\mathbb{ S}^T\balpha\Big\}, \label{comp:alpha}\\
&	\bfgamma^{t+1}=\argmin_{\bfgamma \in {\cal R}^p}\Big\{ R_n(\balpha^{t+1}, \bfgamma)+ \lambda\|\bgamma\|_1\Big\} \label{comp:gamma},
\end{align}
where $R_n(\balpha^{t+1}, \bfgamma):=
\frac{2}{n}({\balpha}^{t+1})^T\mathbb{ S}  \mathbb{K}^c  \mathbb{Z}\bgamma+ {\frac{1}{n}(\by-\mathbb{Z}\bgamma)^T(\by-\mathbb{Z}\bgamma)}$. 

To update $\balpha$,  it is clear that the optimization task \eqref{comp:alpha} has an analytic solution that 
$$
\balpha^{t+1}=\big ( (\mathbb{ S}\mathbb{K}^c)(\mathbb{ S}\mathbb{K}^c)^T+  n\mu^2\mathbb{ S}\mathbb{K}^c\mathbb{ S}'     \big )^{-1}    \mathbb{ S}\mathbb{K}^c    (\by-\mathbb{Z}\bgamma^t) .
$$
To update $\bgamma$,  we first  introduce the proximal operator (Moreau, 1962), which is  defined as
\begin{align}\label{eqn:prox}
\mbox{Prox}_{{\lambda}\|\cdot\|_1}({\bf u}):=\argmin_{{\bf u}} \Big \{  \frac{1}{2}\|{\bf u}-{\bf v}\|_2^2 + \lambda \|{\bf u}\|_1  \Big \}.
\end{align}
Note that the solution of optimization task \eqref{eqn:prox} is the well-known soft-thresholding operator with solution that 
$$
{  \big ( \mbox{Prox}_{{\lambda}\|\cdot\|_1}(\bu) \big )_i=\sign(u_i)(|u_i|-{\lambda})_+}.
$$

Then, for the optimization task \eqref{comp:gamma}, we have
$$
\bfgamma^{t+1}=\mbox{Prox}_{\frac{\lambda}{D}\|\cdot\|_1}\Big (\bfgamma^t -   \frac{1}{D}    \nabla_{\bfgamma} R_n(\balpha^{t+1}, \bfgamma^t)
\Big ),
$$
where   $D$ denotes an upper bound of the  Lipschitz constant of $R_n(\balpha^{t+1}, \bfgamma^t)$,  and compute $\nabla_{\bgamma} R_n(\balpha^{t+1}, \bfgamma^t)= 
\frac{2}{n}  \mathbb{Z}^T (\mathbb{ S}  \mathbb{K}^c )^T {\balpha}^{t+1}    + \frac{2}{n}\mathbb{Z}^T\mathbb{Z}\bgamma^t         - \frac{2}{n}\mathbb{Z}^T\by            $.  
We repeat the above iteration steps  until $(\balpha^{t+1},\bgamma^{t+1})$ converges.

It should be pointed out that the exact value of $ D$ is often difficult to determine in large-scale problems. A common way to handle this problem is to use  a backtracking scheme \citep*{Boyd2004} as a more efficient alternative to approximately compute  an upper bound of it.

\subsection{Choice of Random Sketch Matrix }\label{sec:4.2}
In this paper, we consider  three random sketch methods, including the sub-Gaussian random sketch (GRS), randomized orthogonal system  sketch (ROS) and sub-sampling random sketch (SUB).   Precisely,  we denote the $i$-th row of  the random matrix $\mathbb{S}$ as ${\bs}_i$  and consider three different types of  ${\bs}_i$ as follows.

\noindent {\bf Sub-Gaussian sketch (GRS):}  The row ${\bs}_i$ of  $\mathbb{S}$  is zero-mean $1$-sub-Gaussian if for any $\bu \in {\cal R}^n$, we have
$$
\text{P}\big( \langle {\bs}_i, \bu \rangle  \geq t\|\bu\|_2\big ) \leq e^{-t^2/2}, \ \ \forall \, t \geq 0.
$$
Note that  the row ${\bs}_i$ with independent and identical distributed $N(0,1)$ entries is  1-sub-Gaussian.  For simplicity, we  further rescale the sub-Gaussian sketch matrix $\mathbb{S}$ such that the rows $\bs_i$ have the covariance matrix $\frac{1}{\sqrt{m}}\mathbb{I}_{n}$, where $\mathbb{I}_n$ denotes a $n$ dimensional identity matrix.

\noindent {\bf Randomized orthogonal system  sketch (ROS):}  The row $\bs_i$ of the random matrix $\mathbb{S}$ is formed with i.i.d rows of the form
$$
{\bs}_i=\sqrt{\frac{n}{m}} \mathbb{R}\mathbb{H}^T\mathbb{I}_{(i)},  ~\text{ for }~ i=1,...,m,
$$
where $\mathbb{R}\in {\cal R}^{n\times n}$ is a random diagonal matrix whose entries are i.i.d. Rademacher variables taking value $\{-1,1\}$ with equal probability,  $\mathbb{H}=\{H_{ij}\}_{i,j=1}^n \in {\cal R}^{n\times n}$ is an orthonormal matrix with bounded entries that $H_{ij}\in [-\frac{1}{\sqrt{n}}, \frac{1}{\sqrt{n}}]$, and the $n$-dimensional vectors $\mathbb{I}_{(1)},...,\mathbb{I}_{(m)}$ are drawn uniformly at random without replacement from the $n$-dimensional identity matrix $\mathbb{I}_n$ .

\noindent {\bf Sub-sampling sketches (SUB):}  The rows $\bs_i$ of the random matrix
$\mathbb{S}$ has the form that
$$
{\bs}_i=\sqrt{\frac{n}{m}}\mathbb{I}_{(i)},
$$
where  the $n$-dimensional vectors $\mathbb{I}_{(1)},...,\mathbb{I}_{(m)}$ are drawn uniformly at random without replacement from a $n$ dimensional identity matrix. Note that the  sub-sampling sketches method  can be regarded as a special case of the ROS sketch by replacing the matrix $\mathbb{R}^T\mathbb{H}$ with a $n$-dimensional identity matrix $\mathbb{I}_n$.

\subsection{Choice of the Sketch Dimension}
In practice, we are interested in the $m\times n$ sketch matrices with $m \ll n$ to enhance computational efficiency. Note that the existence of a $n \times n$ kernel matrix in Lemma 1 is only a sufficient condition for equivalent optimization. It has been shown theoretically in the  kernel regression \citep{Yang2017} that the kernel matrix can be compressed to be the one with small size, based on some intrinsic low-dimensional notations. 
Despite  the model difference from \cite{Yang2017}, our kernel matrix  $\mathbb{K}^c$ does not depend on the scalar covariates $\bf Z$, and thus those derived results for the kernel regression are still applicable to our case. 

Consider the eigen-decomposition $\mathbb{K}^c=\mathbb{U}\mathbb{D}\mathbb{U}^T$ of the kernel matrix, where $\mathbb{U}\in {\cal R}^{n\times n}$
is an orthonormal matrix of eigenvectors, and $\mathbb{D}=\hbox{diag}\{\hat{\mu}_1,...,\hat{\mu}_n\}$  is a diagonal
matrix of eigenvalues, where $\hat{\mu}_1\geq \hat{\mu}_2\geq ...\geq\hat{\mu}_n\geq 0$. We define the kernel complexity function as
$$
\widehat{\mathcal{R}}(\delta)=\sqrt{\frac{1}{n}\sum_{j=1}^n\min\{\delta,\hat \mu_j\}}.
$$ 
The critical radius is defined to be the smallest positive solution of $\delta_n>0$ to the inequality
$$
\widehat{\mathcal{R}}(\delta)\leq \delta^2/\sigma.
$$
Note that the existence  and uniqueness of this critical radius is guaranteed for any kernel class. Based on this, we define the statistical dimension of the kernel is
$$
d_n:=\min\{j\in[n]:\hat{\mu}_j\leq \delta^2_n\}.
$$
Recall that, Theorem 2  in \cite{Yang2017} shows that various forms of randomized sketches can achieve the minimax rate using a sketch dimension proportional to the statistical dimension $d_n$. In particular, for Gaussian sketches and ROS sketches, the sketch dimension $m$ is required
satisfy a lower bound of the form 
\[ m\geq 
\begin{cases}
cd_n      &  \hbox{for Gaussian sketches},\\
cd_n\log^4(n) & \hbox{for ROS sketches}.
\end{cases}
\]
Here $c$ is some constant. 
In this paper, we adopt  this specified sketch dimension $m$ to implement our experiments. 


\section{Numerical Experiments}\label{sec:num}

In this section, we illustrate the numerical performance of the proposed method with random sketches in two numerical examples. Specifically, we assume that the true generating model  is
\begin{align}
Y_i=\int_{\cal T} f^*(t)X_i(t)dt+{\bf Z}_i^T\bfgamma^*+\varepsilon_i,
\end{align}
where $\varepsilon_i \sim N(0, \sigma^2)$ with $\sigma=1$, and ${\cal T}$ is set as $[0,1]$. Note that the generating scheme is the same as that in \citealt*{Hall2007} and \citealt*{Yuan2010}. In practice, the integrals in calculation of $\mathbb{B}$ and $\mathbb{K}^c$ are approximated by summations, and thus we generate 1000 points in ${\cal T}=[0,1]$ with equal distance and evaluate the integral by using the generated points. As the proper choice of  tuning parameters  plays a crucial role in achieving the desired performance of the proposed method,  we   adopt 5-fold  cross-validation to select the optimal values of the tuning parameters $\mu$ and $\lambda$. 

In {all} the simulated cases, we consider a RKHS ${\cal H}_K$ induced by a reproducing kernel function on ${\cal T}\times {\cal T}$ that
\begin{align*}
K(s,t) &= \sum_{k\geq 1} \frac{2}{(k\pi)^4} \cos(k\pi s) \cos(k\pi t)\\
&=\sum_{k\geq 1} \frac{1}{(k\pi)^4} \cos(k \pi (s-t)) + \sum_{k\geq 1} \frac{1}{(k\pi)^4} \cos(k \pi (s+t))  \\
&=    -\frac{1}{3} B_4 \big (  \frac{|s-t|}{2} \big ) -\frac{1}{3} B_4 \big (  \frac{s+t}{2}  \big ),
\end{align*}
where $B_{2m}(\cdot)$ denotes the $2m$-th  Bernoulli polynomial that
$$
B_{2m}(s)=  (-1)^{m-1} 2 (2m)! \sum_{k \geq 1} \frac{  \cos(2\pi k s)     }{(2\pi k)^{2m}}, ~\text{for any}~ s\in {\cal T}.
$$
Note that the  RKHS ${\cal H}_K$ induced by  $K(s,t)$ contains the functions in a linear span of the cosine basis that
$$
f(s)=\sqrt{2} \sum_{k\geq 1} g_k \cos(k\pi s), ~\text{for any}~ s\in {\cal T}.
$$
such that $\sum_{k\geq 1} k^4 g_k^2 <\infty$  and the endowed norm is 
$$
\|f\|^2_K=\int_{{\cal T}} \big ( \sqrt{2} \sum_{k\geq 1} (k\pi)^2g_k\cos(k\pi t) \big )^2 dt =\sum_{k\geq 1} (k\pi)^4g_k^2.
$$


The performance of the proposed method is evaluated under the following two numerical examples.

\noindent {\bf Example 1}. We consider  the true slope function  $f^*$ and the random function $X$ are
$$
f^*(t)=\sum_{k=1}^{50}4(-1)^{k+1}k^{-2}\sqrt{2}\cos(k\pi t),
$$ 
and  
$$
X(t)=\xi_1U_1+\sum_{k=2}^{50} \xi_kU_k \sqrt{2} \cos (k \pi t),
$$
where $U_k \sim U(-\sqrt{3}, \sqrt{3})$ and   $\xi_k=(-1)^{k+1}k^{-v/2}$. For the linear part, the true regression coefficients are set as $\bfgamma^0=(2,-2,0,...,0)^T$ and the sample $\mathbb{ Z}=({\bf Z}_1,...,{\bf Z}_n)^T\in {\cal R}^{n\times p}$ with ${\bf Z}_i=(z_{i1},...,z_{ip})^T$ are generated i.i.d. as $z_{ij} \sim U(0,1)$. 

\noindent{\bf Example 2.} The generating scheme is the same as  Example 1, except that
\[ \xi_k=
\begin{cases}
1,      & \quad k=1,\\
0.2(-1)^{k+1}(1-0.0001k),  & \quad 2\leq k \leq 4,\\
0.2(-1)^{k+1}\big [   (5\lfloor  k/5\rfloor)^{-\upsilon/2}- 0.0001(k ~\text{mod}~ 5) \big ],  & \quad k \geq  5.
\end{cases}
\]

Clearly, $\xi^2_k$'s are the  eigenvalues of the covariance function   $C$ and we  choose $v=1.1, 2$ and $4$ to evaluate the effect of the smoothness of $\xi_k$ in the both examples. Note that in Example 1, these eigenvalues are well spaced, and  the covariance function   $C$ and the reproducing kernel $K$ share the same eigenvalues, while in Example 2, these eigenvalues are closely spaced and the alignment between $K$ and $C$ is considered.

\begin{figure}[!h]
	\centering
	\begin{subfigure}[b]{0.3\textwidth}
		\includegraphics[width=\textwidth]{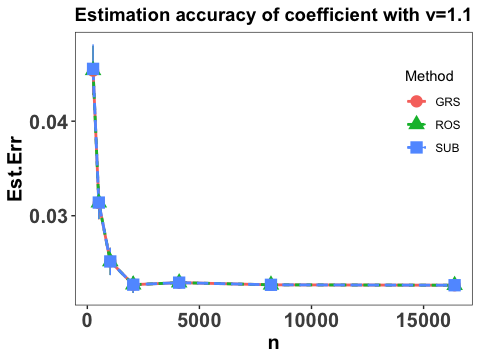}
	\end{subfigure}
	\begin{subfigure}[b]{0.3\textwidth}
		\includegraphics[width=\textwidth]{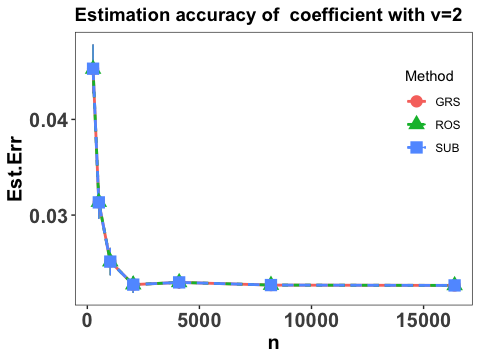}
	\end{subfigure}
	\begin{subfigure}[b]{0.3\textwidth}
		\includegraphics[width=\textwidth]{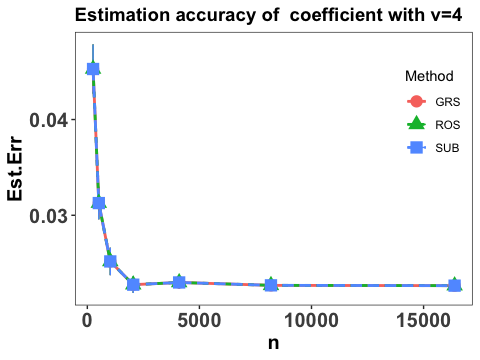}
	\end{subfigure}
	\caption{Estimation accuracy of the coefficients in Example 1 under various scenarios.	}
	\label{fig:101}
\end{figure}

To comprehend the effect of sample size, we 
consider the same settings as in \cite{Yang2017} that  $n=256, 512, 1024, 2048, 4096, 8192$ and $16384$ and conservatively, take $m= \lfloor n^{1/3} \rfloor$ for the three random sketch methods introduced in Section \ref{sec:4.2}. Note that with the choice of  $m$, the time and store complexities reduce to  $O(n)$ and $O(n^{2/3})$, respectively. Each scenario is replicated 50 times and the performance of the proposed method is evaluated by various measures, including the estimation accuracy of the linear coefficients, the integrated prediction error in terms of the slope function and the  prediction error.  Specifically,   the estimation accuracy of the linear coefficients   is evaluated by 
$
\|\widehat{\bfgamma}-\bfgamma^0\|^2_2=\sum_{l=1}^p (\widehat{\gamma}_l-\gamma_l^0)^2,
$
and Figure  \ref{fig:101} shows the estimation accuracy of the coefficients with different choice of $v$.

It is clear that the estimation error of the coefficients converges linearly as sample size $n$ increases and becomes stable when $n$ is sufficiently large, and  the three employed sketch methods have similar performance. It is also interesting to notice that  the convergence patterns under  difference choice of $v$ are almost the same,  which concurs with our theoretical findings that estimation of $\bgamma^*$ is not affected by the existence of nonparametric components in Theorems \ref{nonpara} and \ref{thm3}. 

\begin{figure}[!h]
	\centering
	\begin{subfigure}[b]{0.3\textwidth}
		\includegraphics[width=\textwidth]{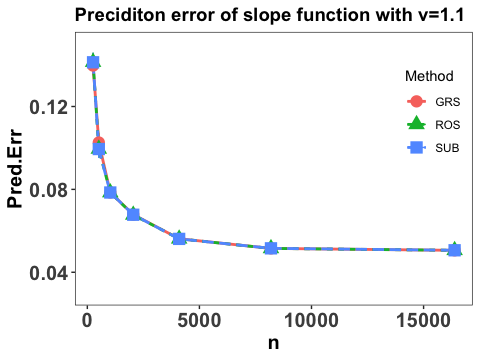}
	\end{subfigure}
	\begin{subfigure}[b]{0.3\textwidth}
		\includegraphics[width=\textwidth]{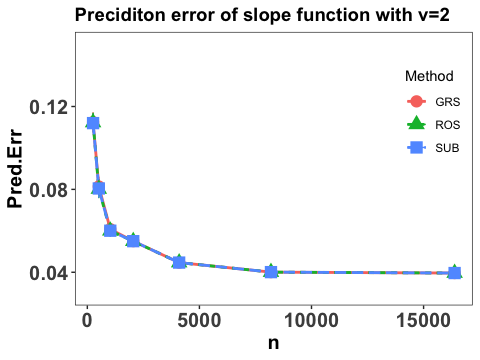}
	\end{subfigure}
	\begin{subfigure}[b]{0.3\textwidth}
		\includegraphics[width=\textwidth]{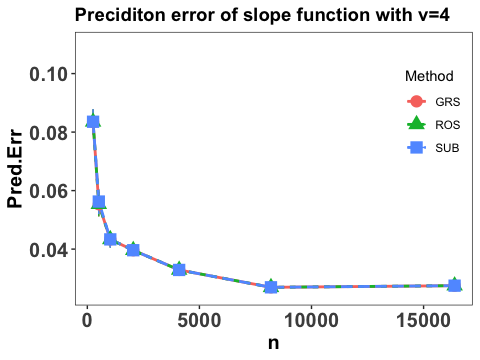}
	\end{subfigure}
	\caption{Prediction error of the slope function in Example 1 under various scenarios.	}
	\label{fig:102}
\end{figure}

Let $(Y', {X'}(\cdot), {\bf Z}')$ denotes an independent copy of $(Y, {X}(\cdot), {\bf Z})$ and the integrated prediction error in terms of the slope function is reported by
$$
\widehat{\mathbb{E}}_{X'}\|\widehat{f}-f^*\|^2= \widehat{\mathbb{E}}_{X'}   \big (  \int_{\cal T} ( \widehat{f}(t)-f^*(t) )X'(t) dt \big  )^2
$$ 
The empirical expectation  $\widehat{\mathbb{E}}$ is evaluated by  a testing sample with size $10000$ and   $\widehat{Y}'=  \int_{\cal T} \widehat{f}(t)X^{'}_i(t)dt+({\bf Z}^{'}_i)^T\widehat{\bfgamma}$ and the numerical performance  are summarized in Figure \ref{fig:102}. 

Note that Figure  \ref{fig:102} suggests  that the prediction error of the slope function converges at some polynomial rate as sample size $n$, which agrees with our theoretical results in Section \ref{sec:3}, and the three employed sketch methods yield similar numerical performance. Moreover,  it can be seen that with the increase of the value of $v$, the prediction error  goes down, which also concurs with our theoretical findings in Theorems \ref{upperb} and \ref{thm3} that the faster decay rate of the eigenvalues, the smaller the prediction error.

\begin{figure}[!h]
	\centering
	\begin{subfigure}[b]{0.3\textwidth}
		\includegraphics[width=\textwidth]{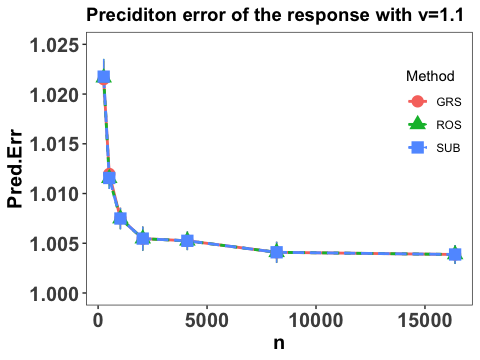}
	\end{subfigure}
	\begin{subfigure}[b]{0.3\textwidth}
		\includegraphics[width=\textwidth]{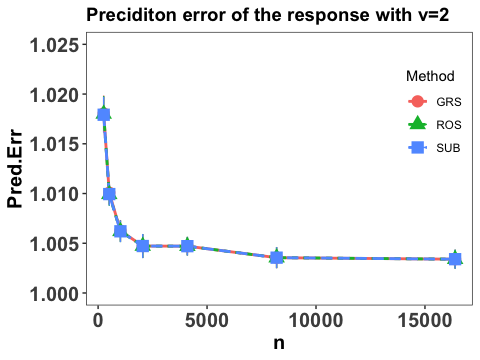}
	\end{subfigure}
	\begin{subfigure}[b]{0.3\textwidth}
		\includegraphics[width=\textwidth]{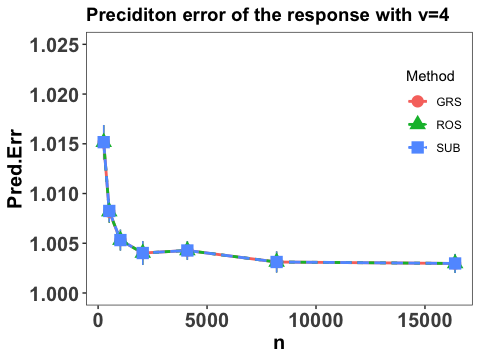}
	\end{subfigure}
	\caption{Prediction error of the response in Example 1 under various scenarios.	}
	\label{fig:103}
\end{figure}

We also report  the integrated    prediction error of the response  by calculating
$$
\widehat{\mathbb{E}}_{Y', X'}\|\widehat{Y}'-Y'\|^2_2.
$$
The empirical expectation  $\widehat{\mathbb{E}}$ is also evaluated by  a testing sample with size $10000$  and the numerical performance  are summarized in Figure \ref{fig:103}.

Clearly, we conclude that prediction error of the response converges at some polynomial rate as sample size $n$ and the prediction error becomes smaller with $v$ increases, which agrees with our theoretical results in Theorem \ref{upperb}. It is also interesting to point out that  the three employed sketch methods yield similar numerical performance and the prediction errors tends to converge to 1, which is the variance of $\varepsilon$ in the true modelling. This verifies  the efficiency of the proposed estimation and the proper choice of $m$.   

\begin{figure}[!h]
	\centering
	\begin{subfigure}[b]{0.3\textwidth}
		\includegraphics[width=\textwidth]{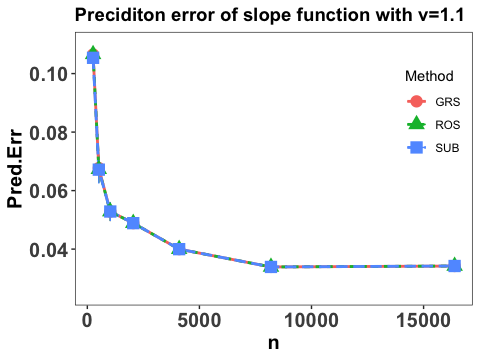}
	\end{subfigure}
	\begin{subfigure}[b]{0.3\textwidth}
		\includegraphics[width=\textwidth]{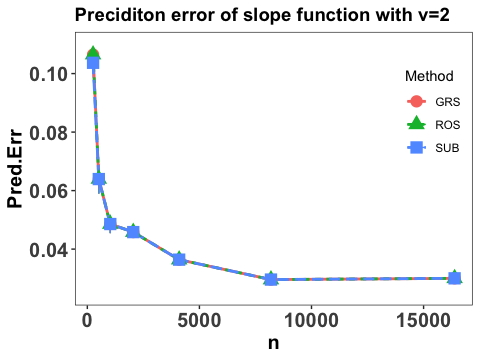}
	\end{subfigure}
	\begin{subfigure}[b]{0.3\textwidth}
		\includegraphics[width=\textwidth]{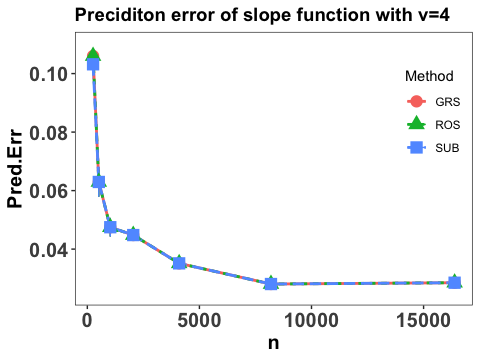}
	\end{subfigure}
	\begin{subfigure}[b]{0.3\textwidth}
		\includegraphics[width=\textwidth]{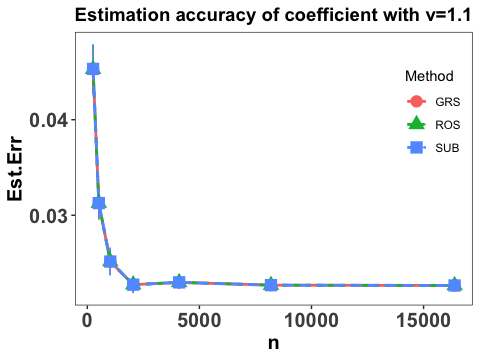}
	\end{subfigure}
	\begin{subfigure}[b]{0.3\textwidth}
		\includegraphics[width=\textwidth]{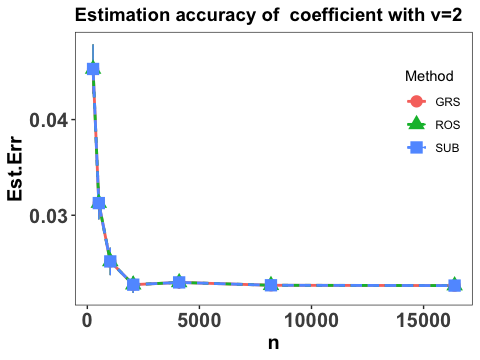}
	\end{subfigure}
	\begin{subfigure}[b]{0.3\textwidth}
		\includegraphics[width=\textwidth]{EX23co.png}
	\end{subfigure}
	\begin{subfigure}[b]{0.3\textwidth}
		\includegraphics[width=\textwidth]{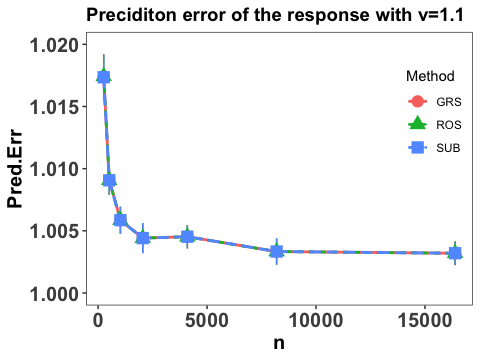}
	\end{subfigure}
	\begin{subfigure}[b]{0.3\textwidth}
		\includegraphics[width=\textwidth]{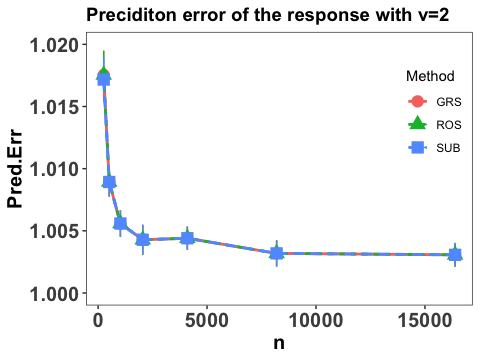}
	\end{subfigure}
	\begin{subfigure}[b]{0.3\textwidth}
		\includegraphics[width=\textwidth]{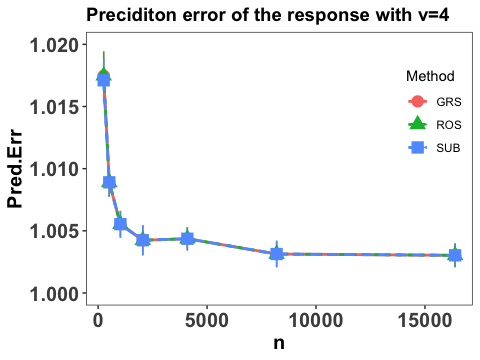}
	\end{subfigure}
	\caption{Numerical performance of the proposed method in Example 2 under various scenarios.	}
	\label{fig:21}
\end{figure}

Note that the numerical results in Example 2 where the eigenvalues are closely spaced are  similar to those in the case with well-spaced eigenvalues in Example 1. 
Figure \ref{fig:21} shows the numerical performance under  the closely spaced eigenvalues setting in Example 2. 


\section{Conclusion}\label{sec:con}
This paper establishes the optimal minimax rate for the estimation of partially functional linear model (PFLM) under kernel-based and high dimensional setting. The optimal minimax rates of estimation is established by using various techniques in empirical process theory for analyzing kernel classes, and an efficient numerical algorithm based on randomized sketches of the kernel matrix is implemented to verify our theoretical findings.

\acks{Shaogao Lv's research was partially supported by NSFC-11871277.  Xin He’s research was supported in part by NSFC-11901375 and Shanghai Pujiang Program 2019PJC051.  Junhui Wang’s research was supported in part by HK RGC Grants GRF-11303918  and GRF-11300919.}


\bibliography{bibtex.bib}

\begin{thebibliography}{25}
\providecommand{\natexlab}[1]{#1}
\providecommand{\url}[1]{\texttt{#1}}
\expandafter\ifx\csname urlstyle\endcsname\relax
  \providecommand{\doi}[1]{doi: #1}\else
  \providecommand{\doi}{doi: \begingroup \urlstyle{rm}\Url}\fi

\bibitem[Aronszajn(1950)]{Aronszajn1950}
N.~Aronszajn.
\newblock Theory of reporudcing kernels.
\newblock \emph{Transactions of the American Mathematical Society},
  \textbf{68}:\penalty0 337--404, 1950.

\bibitem[Bickel et~al.(2009)Bickel, Ritov, and Tsybakov]{Bickel2009}
P.~Bickel, Y.~Ritov, and A.~Tsybakov.
\newblock Simultaneous analysis of lasso and dantzig selector.
\newblock \emph{Annals of Statistics}, \textbf{37}:\penalty0 1705--1732, 2009.

\bibitem[Bousquet(2002)]{Bousquet2002}
O.~Bousquet.
\newblock A bennet concentration inequality and its application to suprema of
  empirical processes.
\newblock \emph{Comptes Rendus Mathematique Academie des Sciences Paris},
  \textbf{334}:\penalty0 495--550, 2002.

\bibitem[Boyd and Vandenberghe(2004)]{Boyd2004}
S.~Boyd and L.~Vandenberghe.
\newblock \emph{Convex Optimization}.
\newblock Cambridge University Press, Cambridge, 2004.

\bibitem[B\"uhlmann and Van.~de. Geer(2011)]{Buhlmann2011}
P.~B\"uhlmann and S.~Van.~de. Geer.
\newblock \emph{Statistics for High-Dimensional Data: Methods, Theory and
  Applications}.
\newblock Springer, Heidelberg, 2011.

\bibitem[Cai and Yuan(2012)]{Tony2012}
T.~Cai and M.~Yuan.
\newblock Minimax and adaptive prediction for functional linear regression.
\newblock \emph{Journal of the American Statistical Association},
  \textbf{107}:\penalty0 1201--1216, 2012.

\bibitem[Cardot et~al.(2003)Cardot, Ferraty, and Sarda]{Cardot2003}
H.~Cardot, F.~Ferraty, and P.~Sarda.
\newblock Spline estimators for the functional linear model.
\newblock \emph{Statistica Sinica}, \textbf{13}:\penalty0 571--591, 2003.

\bibitem[Ferraty and Vieu(2006)]{Ferraty2006}
F.~Ferraty and P.~Vieu.
\newblock \emph{Nonparametric Functional Data Analysis: Theory and Practice}.
\newblock Springer, New York, 2006.

\bibitem[Hall and Horowitz(2007)]{Hall2007}
P.~Hall and J.~Horowitz.
\newblock Methodology and convergence rates for functional linear regression.
\newblock \emph{Annals of Statistics}, \textbf{35}:\penalty0 70--91, 2007.

\bibitem[Kong et~al.(2016)Kong, Xue, Yao, and Zhang]{Kong2016}
D.~Kong, K.~Xue, F.~Yao, and H.~Zhang.
\newblock Partially functional linear regression in high dimensions.
\newblock \emph{Biometrika}, \textbf{103}:\penalty0 1--13, 2016.

\bibitem[Ledoux(1997)]{Ledoux1997}
M.~Ledoux.
\newblock On talagrand's deviation inequalities for product measures.
\newblock \emph{Probability and Statistics}, \textbf{1}:\penalty0 63--87, 1997.

\bibitem[Ledoux(2001)]{Ledoux2001}
M.~Ledoux.
\newblock \emph{The Concentration of Measure Phenomenon (Mathematical Surveys
  and Monographs)}.
\newblock American Mathematical Society, Providence, RI, 2001.

\bibitem[Lu et~al.(2014)Lu, Du, and Sun]{Lu2014}
Y.~Lu, J.~Du, and Z.~Sun.
\newblock Functional partially linear quantile regression model.
\newblock \emph{Metrika}, \textbf{77}:\penalty0 17--32, 2014.

\bibitem[Mahoney(2011)]{Mahoney2011}
M.~Mahoney.
\newblock Randomized algorithms for matrices and data.
\newblock \emph{Foundations and Trends in Machine Learning in Machine
  Learning}, \textbf{3}:\penalty0 1--54, 2011.

\bibitem[M\"uller and Van~de Geer(2015)]{Muller2015}
P.~M\"uller and S.~Van~de Geer.
\newblock The partial linear model in high dimensions.
\newblock \emph{Scandinavian Journal of Statistics}, \textbf{42}:\penalty0
  580--608, 2015.

\bibitem[Ramsay and Silverman(2005)]{Ramsay2005}
J.~Ramsay and B.~Silverman.
\newblock \emph{Functional Data Analysis}.
\newblock Springer, New York, 2005.

\bibitem[Shin and Lee(2012)]{Shin2012}
H.~Shin and M.~Lee.
\newblock On prediction rate in partial functional linear regression.
\newblock \emph{Journal of Multivariate Analysis}, \textbf{103}:\penalty0
  93--106, 2012.

\bibitem[Tsybakov(2009)]{Tsybakov2009}
A.~Tsybakov.
\newblock \emph{Introduction to Nonparametric Estimation}.
\newblock Springer, New York, 2009.

\bibitem[Van.~de. Geer(2000)]{Van2000}
S.~Van.~de. Geer.
\newblock \emph{Emprical Processes in M-Estimation}.
\newblock Cambridge University Press, New York, 2000.

\bibitem[Verzelen(2012)]{Verzelen2012}
N.~Verzelen.
\newblock Minimax risks for sparse regressions: Ultra-high dimensional
  phenomenons.
\newblock \emph{Electronic Journal of Statistics}, \textbf{6}:\penalty0
  38--–90, 2012.

\bibitem[Yang and Barron(1999)]{Yang1999}
Y.~Yang and A.~Barron.
\newblock Information-theoretic determination of minimax rates of convergence.
\newblock \emph{Annals of Statistics}, \textbf{27}:\penalty0 1564--1599, 1999.

\bibitem[Yang et~al.(2017)Yang, Pilanci, and Wainwright]{Yang2017}
Y.~Yang, M.~Pilanci, and M.~Wainwright.
\newblock Randomized sketches for kernels: fast and optimal non-parametric
  regression.
\newblock \emph{Annals of Statistics}, \textbf{45}:\penalty0 991--1023, 2017.

\bibitem[Yu et~al.(2019)Yu, Levine, and Cheng]{Yu2018}
Z.~Yu, M.~Levine, and G.~Cheng.
\newblock Minimax optimal estimation in partially linear additive models under
  high dimension.
\newblock \emph{Bernoulli}, \textbf{25}:\penalty0 1289--1325, 2019.

\bibitem[Yuan and Cai(2010)]{Yuan2010}
M.~Yuan and T.~Cai.
\newblock A reproducing kernel hilbert space approach to functional linear
  regression.
\newblock \emph{Annals of Statistics}, \textbf{38}:\penalty0 3412--3444, 2010.

\bibitem[Zhu et~al.(2014)Zhu, Yao, and Zhang]{Zhu2014}
H.~Zhu, F.~Yao, and H.~Zhang.
\newblock Structured functional additive regression in reproducing kernel
  hilbert spaces.
\newblock \emph{Journal of the Royal Statistical Society, Series B},
  \textbf{76}:\penalty0 581--603, 2014.

\end{thebibliography}

\end{document}